\documentclass[12pt]{article}
\usepackage[english]{babel}
\usepackage{amsmath}
\usepackage[cp1251]{inputenc}
\usepackage{amscd}
\usepackage{amssymb}
\newtheorem{thm}{Theorem}
\newtheorem{lem}{Lemma}

\newtheorem{cor}{Corollary}

\begin{document}

\begin{center}
\bf On reverse Markov-Nikol'skii inequalities\\
for polynomials with restricted zeros
\end{center}

\begin{center}
{\bf Mikhail A. Komarov}\footnote{Vladimir State University, Gor'kogo street 87,
600000 Vladimir, Russia; \ \ e-mail: kami9@yandex.ru}
\end{center}

{\bf Abstract.} {\small Let $\Pi_n$ be the class of
algebraic polynomials $P$ of degree $n$,
all of whose zeros lie on the segment $[-1,1]$.
In 1995, S.~P.~Zhou has proved the following Tur\'{a}n type
reverse Markov-Nikol'skii inequality:
$\|P'\|_{L_p[-1,1]}>c\, {(\sqrt{n})}^{1-1/p+1/q}\, \|P\|_{L_q[-1,1]}$, $P\in \Pi_n$,
where $0<p\le q\le \infty$, $1-1/p+1/q\ge 0$
($c>0$ is a constant independent of $P$ and $n$).
We show that Zhou's estimate remains true in the case
$p=\infty$, $q>1$. Some of related Tur\'{a}n
type inequalities are also discussed.\\
}

For any $n\in \mathbb{N}$, let $\Pi_n$ be the class of polynomials $P$ of degree $n$,
all of whose zeros lie in the closed interval $[-1,1]$,
and $\Pi_n(D^+)$ --- the class of polynomials $P$ of degree $n$,
all of whose zeros lie in the closed half-disk
$D^+:=\{z:|z|\le 1,\ {\rm Im}\,z\ge 0\}$. We will use
the notation $\|P\|_p:=\|P\|_{L_p[-1,1]}$,
$0<p\le \infty$, where $\|P\|_{L_\infty[-1,1]}=\|P\|_{C[-1,1]}$.

In \cite{K-AnalysisMath}, the author proved the inequality
\begin{equation}\label{result-AnalysisMath}
    \|P'\|_\infty>A\sqrt{n}\,\, \|P\|_\infty, \quad
    A=\frac{2}{3\sqrt{210 e}}, \quad P\in \Pi_n(D^+),
\end{equation}
that generalizes the well-known {\it reverse Markov inequality}
\begin{equation}\label{Turan}
    \|P'\|_\infty>\frac{1}{6}\sqrt{n}\,\, \|P\|_\infty, \quad
    P\in \Pi_n,
\end{equation}
obtained by P.~Tur\'{a}n in 1939. The $L_p$-version of (\ref{Turan})
was obtained for $p=2$ --- by A.~K.~Varma (1976) and for any $p>0$ --- by
S.~P.~Zhou (1984, 1986), see \cite{Varma76}, \cite{Zhou86}:
\begin{equation}\label{Zhou86}
    \|P'\|_p>c\,\sqrt{n}\,\, \|P\|_p, \ \ P\in \Pi_n,
    \quad 0<p<\infty,
\end{equation}
where $c>0$ is a constant depending only on $p$.
The following more general result ({\it reverse Markov-Nikol'skii inequality})
is also due to Zhou \cite{Zhou95}: suppose that
\[0<p\le q\le \infty, \quad 1-1/p+1/q\ge 0;\]
then
\begin{equation}\label{Zhou95}
    \|P'\|_p>c\, {(\sqrt{n})}^{1-1/p+1/q}\, \|P\|_q, \quad P\in \Pi_n,
\end{equation}
where again $c=c(p)>0$.
Taking the polynomials $Q_n(x):=(1-x^2)^n$, we see that the order of $n$
in (\ref{Zhou95}) (and hence in (\ref{result-AnalysisMath})-(\ref{Zhou86}))
can not be improved, because
\[\int_{-1}^1 {|Q_n'(x)|}^p\, dx\asymp n^{(p-1)/2},
\quad \int_{-1}^1 {|Q_n(x)|}^q\, dx\asymp n^{-1/2}
\quad (n\to \infty)\]
for any fixed $p,q>0$.

For extensions of (\ref{result-AnalysisMath}) and (\ref{Zhou95}),
see \cite{Erdelyi2020} and \cite{Zhao-Zhou-Yu-Wang}, correspondingly.
The limiting case $p=0$ of (\ref{Zhou86}) was investigated in \cite{K-AnalysisMath21}.
We also recall the following weighted Tur\'{a}n type inequalities
with the best possible constants: if $P\in \Pi_n$, then
\begin{equation}\label{BabenkoPichugov}
    \big\|P'(x)\sqrt{1-x^2}\big\|_\infty\ge A_\infty(n)\, \|P\|_\infty, \quad
    A_\infty(n):=\frac{\sqrt{n}}{\sqrt{2}}{\Big(1-\frac{1}{2n}\Big)}^{n-1/2}
\end{equation}
(V.~F.~Babenko and S.~A.~Pichugov, \cite{BabenkoPichugov(MathNotes)}), and
\begin{equation}\label{Tyrygin-UMJ}
    \big\|P'(x)(1-x^2)^\frac{p-1}{2p}\big\|_p\ge A_p(n)\, \|P\|_\infty,
    \quad 1\le p<\infty,
\end{equation}
\[A_p(n):=n\, \Big\{\mathrm{B}\left(\frac{(2n-1)p+1}{2},\frac{p+1}{2}\right)\Big\}^{1/p}\]
(I.~Ya.~Tyrygin, \cite{Tyrygin-UMJ}), where $\mathrm{B}(x,y)$
is the Euler beta-function. It can be checked that
$A_p(n)\to A_\infty(n)$ ($p\to \infty$) and
\[A_p(n)\asymp {(\sqrt{n})}^{1-1/p}, \quad 1\le p\le \infty
\quad (n\to \infty).\]

In this short note, we show the existence of the estimate of the form
(\ref{Zhou95}) in the case when $p=\infty$, $q>1$ (in contrast
to Zhou's case, where $p\le q$). Namely, for the class $\Pi_n(D^+)$,
the following extension of (\ref{Zhou95}), as well as of (\ref{result-AnalysisMath}), holds:

\begin{thm}\label{Theorem 1}
If $P\in \Pi_n(D^+)$, then for any $1<q\le \infty$
\begin{equation}\label{Theorem D+}
    \|P'\|_\infty>6^{1/q}(1-q^{-1})^{1/q}\, {(A\sqrt{n})}^{1+1/q}\, \|P\|_q,
\end{equation}
where $A$ is the constant defined in $(\ref{result-AnalysisMath})$.
\end{thm}

This bound can be slightly improved for the class $\Pi_n$ as follows:

\begin{thm}\label{Theorem 2}
If $P\in \Pi_n$, then for any $1<q\le \infty$
\begin{equation}\label{Theorem [-1,1]}
    \big\|P'(x)\sqrt{1-x^2}\big\|_\infty^{1-1/q}\,
    \big\|P'(x)(1-x^2)\big\|_\infty^{1/q}
    \ge c_q\,{n}^{1/q}\, \{A_\infty(n)\}^{1-1/q}\,
    \|P\|_q,
\end{equation}
where $A_\infty(n)>\sqrt{n}/\sqrt{2e}$ is the constant defined in $(\ref{BabenkoPichugov})$,
$c_q:={2}^{-1/q}(1-q^{-1})^{1/q}$.
\end{thm}

\begin{cor}\label{Cor 1}
If $P\in \Pi_n$, then for any $1<q\le \infty$
\[\big\|P'(x)\sqrt{1-x^2}\big\|_\infty
> c_q\, {(\sqrt{2e})}^{-1+1/q}\, {(\sqrt{n})}^{1+1/q}\, \|P\|_q.\]
\end{cor}

Of course, in the case $q=\infty$ we have (\ref{Theorem D+}) and (\ref{Theorem [-1,1]})
directly from (\ref{result-AnalysisMath}) and (\ref{BabenkoPichugov}).

Next, applying (\ref{BabenkoPichugov}), (\ref{Tyrygin-UMJ})
to the derivative $P'\in \Pi_{n-1}$
and combining the result with (\ref{Theorem [-1,1]}), we get the reverse
Markov-Nikol'skii inequality for the second derivative.

\begin{cor}\label{Cor 2}
If $P\in \Pi_n$, $n\ge 2$, then for any $1\le p\le \infty$
and $1<q\le \infty$
\[\big\|P''(x)(1-x^2)^\frac{p-1}{2p}\big\|_p\ge C_{p,q}(n)\, \|P\|_q,\]
where $C_{p,q}(n):=c_q\, {n}^{1/q}\, \{A_\infty(n)\}^{1-1/q}\, \{A_p(n-1)\}$,
\[C_{p,q}(n)\asymp {(\sqrt{n})}^{2-1/p+1/q} \quad (n\to \infty).\]
\end{cor}

For example, $C_{1,\infty}=A_\infty$, and we have
$\|P''\|_1>(\sqrt{n}/\sqrt{2e})\, \|P\|_\infty$.

It should be mentioned that the inequality $\|P''\|_\infty\ge C(n)\,\|P\|_\infty$
with the best possible constant $C(n)=\min\{n;(n-1)n/4\}$
was constructed by Babenko and Pichugov \cite{BabenkoPichugov(UkrMathJ)}.

Taking again $Q_n(x)=(1-x^2)^n$, we see that Theorems \ref{Theorem 1}, \ref{Theorem 2}
and Corollaries \ref{Cor 1}, \ref{Cor 2} can not be improved in the order of $n$.
In the case of Corollary \ref{Cor 2} we use the fact that
$Q_n''(x)=2n(1-x^2)^{n-2}\{x^2(2n-1)-1\}$, where
$|x^2(2n-1)-1|\le 1$ ($|x|\le 1/\sqrt{2n-1}$)
and $0<x^2(2n-1)-1<2nx^2$ ($|x|>1/\sqrt{2n-1}$). Hence
\[\int_{-1}^1 {|Q_n''(x)|}^p\, dx< 2\int_0^{1/\sqrt{2n-1}}(2n)^p\,dx
+2\int_{1/\sqrt{2n-1}}^1 {(2nx)}^{2p}{(1-x^2)}^{(n-2)p}\,dx\]
\[< 2^{1+p}\, n^{p-1/2}
+2^{1+2p}\,n^{2p}\int_0^1 x^{2p}{(1-x^2)}^{(n-2)p}\,dx
\asymp n^{p-1/2}.\]

{\it Proof of Theorems \ref{Theorem 1} and \ref{Theorem 2}.}
To prove the theorems, we need the following non-standard
Nikol'skii type inequality between different metrics:

\begin{lem}\label{Lem}
Let $1<q<\infty$. Then
\begin{equation}\label{Lemma}
    \int_{-1}^1 {|Q(x)|}^q\, dx\le K\frac{q}{q-1}\, \frac{{\|Q\|}^{q-1}_\infty}{n}
\end{equation}
with $K:=70e$ for any polynomial $Q$ of degree $n$ such that
\[ Q\in \Pi_n(D^+), \quad \|Q'\|_\infty=1.\]
Moreover, $(\ref{Lemma})$ holds true with $K:=2$, if
\begin{equation}\label{type of Q}
    Q\in \Pi_n, \quad \left\|Q'(x)(1-x^2)\right\|_\infty=1.
\end{equation}
\end{lem}

Thus, from Lemma \ref{Lem} we immediately deduce that for any polynomial $P\in \Pi_n$
\begin{equation}\label{after Lemma}
    n\,2^{-1}\big(1-q^{-1}\big) \int_{-1}^1 {|P(x)|}^q\, dx\le
    \left\|P'(x)(1-x^2)\right\|_\infty {\|P\|}^{q-1}_\infty
    \quad (q>1),
\end{equation}
and (\ref{Theorem [-1,1]}) follows by this and
$\|P\|_\infty\le \left\|P'(x)\sqrt{1-x^2}\right\|_\infty/A_\infty$
(see (\ref{BabenkoPichugov})). The proof of (\ref{Theorem D+})
is similar and uses the inequalities $\|P\|_\infty<\|P'\|_\infty/(A\sqrt{n})$,
$P\in \Pi_n(D^+)$ (see (\ref{result-AnalysisMath})),
and $K=70e=4/(27A^{2})<1/(6 A^{2})$.

\smallskip

{\it Remark.} Some ``quasi'' Tur\'{a}n type estimates
for the class $\Pi_n$ follow from (\ref{after Lemma})
with (\ref{Tyrygin-UMJ}). For example,
we have $c_q\, n^{1/q}\, \{A_p(n)\}^{1-1/q}\, \|P\|_q\le
{\|P'\|}^{1/q}_\infty\, {\|P'\|}^{1-1/q}_p$
$(p\ge 1,\ q>1)$.

\medskip

{\it Proof of Lemma \ref{Lem}.} If a given polynomial $Q$ from $\Pi_n(D^+)$
(or $\Pi_n$) satisfies $n\|Q\|_\infty\le 1$, then (\ref{Lemma}) obviously holds with $K=2(q-1)/q<2$:
\[\int_{-1}^1 {|Q(x)|}^q\, dx\le 2{\|Q\|}^q_\infty\le \frac{2}{n}\,{\|Q\|}^{q-1}_\infty.\]

Now let $Q\in \Pi_n$ be a polynomial such that (\ref{type of Q})
holds and $n\|Q\|_\infty>1$.

For any $\delta>0$, let us define the set
\[E(\delta)=\{x\in [-1,1]: n|Q(x)|>\delta\}.\]
The measure $m(E(\delta))=:h(\delta)$ is zero if $\delta\ge T:=n\|Q\|_\infty$, therefore
\[\int_{-1}^1 \big(n|Q(x)|\big)^q\, dx=
q\int_0^\infty \delta^{q-1} h(\delta)\,d\delta=
q\int_0^T \delta^{q-1} h(\delta)\,d\delta\]
\[\le 2q\int_0^1 \delta^{q-1}\,d\delta +q\int_1^T \delta^{q-1} h(\delta)\,d\delta
=2+q\int_1^T \delta^{q-1} h(\delta)\,d\delta.\]
By (\ref{type of Q}), we have $(1-x^2)|Q'(x)|\le 1$ for $-1\le x\le 1$.
Hence
\begin{equation}\label{E subset}
    E(\delta)\subset \Big\{x\in [-1,1]:
    \frac{(1-x^2)|Q'(x)|}{n |Q(x)|}<\frac{1}{\delta}\Big\}.
\end{equation}
But for any polynomial $\pi_n\in \Pi_n$ and any $\alpha>0$,
the measure
\[m\Big\{x\in [-1,1]: (1-x^2)\left|\frac{\pi_n'(x)}{n\pi_n(x)}\right|\le \alpha\Big\}\le
2\alpha;\]
this was shown in \cite{K-TrudyIMM} (see also \cite{K-mz-2022})
with help of P.~Borwein's identity \cite[p. 75]{Borwein}.
This metric estimate, applied to (\ref{E subset}), provides the bound
\[h(\delta)\le 2/\delta,\]
and we easily get the desired inequality (\ref{Lemma}):
\[n^q \int_{-1}^1 {|Q(x)|}^q\, dx
\le 2+2q\int_1^T \delta^{q-2}\,d\delta
=2+2q\frac{{T}^{q-1}-1}{q-1}<\frac{2q {T}^{q-1}}{q-1}, \quad T=n\|Q\|_\infty.\]

To obtain (\ref{Lemma}) in the remaining case
($Q\in \Pi_n(D^+)$, $\|Q'\|_\infty=1$, $n\|Q\|_\infty>1$), we use
similar arguments with another metric estimate
\[m\Big\{x\in [-1,1]: \left|\frac{\pi_n'(x)}{n \pi_n(x)}\right|\le \alpha\Big\}<K\alpha,
\quad K=70e, \quad \pi_n\in \Pi_n(D^+)\]
(N.~V.~Govorov and Yu.~P.~Lapenko, \cite{Govorov-Lapenko}).
Lemma is proved.

\smallskip

{\it Remark.} Concerning the case $q=1$, it seems that
(for example) the right version of (\ref{Theorem D+})
should be of the form $\|P'\|_\infty>C n\,\|P\|_1$.
Our approach gives the weaker bound
\[\|P'\|_\infty>(30e)^{-1}\,\frac{n}{9+\ln n}\,\|P\|_1, \quad P\in \Pi_n(D^+).\]
Indeed (see (\ref{result-AnalysisMath}) and the proof of Lemma \ref{Lem}), in this case
the polynomial $Q(x):=P(x)/\|P'\|_\infty$
satisfies $\|Q'\|_\infty=1$,
$T=n\|Q\|_\infty<n\, (A\sqrt{n})^{-1}=\sqrt{n}/A$ and
\[n\int_{-1}^1 |Q(x)|\, dx < 2+\int_1^{\sqrt{n}/A}
\frac{70e}{\delta}\,d\delta=2+70e \ln (\sqrt{n}/A)<30e(9+\ln n).\]

\end{document}